\pdfoutput=1
\documentclass[12pt]{article}

\usepackage{amsmath,amssymb,amsthm}
\usepackage{array,longtable}
\usepackage{float}
\usepackage{algorithm}
\usepackage{algpseudocode}
\usepackage[font=small,labelfont=bf]{caption}
\usepackage{booktabs}
\usepackage{hyperref}
\hypersetup{colorlinks=true,linkcolor=blue,citecolor=blue,urlcolor=blue}

\newcommand{\doi}[1]{doi:\,\href{https://doi.org/#1}{\nolinkurl{#1}}}
\newcommand{\TC}{\mathrm{TC}}
\newcommand{\SC}{\mathrm{SC}}
\newcommand{\SCst}{\SC_{\mathrm{strict}}}
\newcommand{\gen}[1]{\langle #1\rangle}
\newcommand{\Zt}{\mathbb{Z}_2}

\newtheorem{theo}{Theorem}[section]
\newtheorem{lemma}[theo]{Lemma}
\newtheorem{corollary}[theo]{Corollary}
\newtheorem{proposition}[theo]{Proposition}
\newtheorem{definition}[theo]{Definition}
\newtheorem{remark}[theo]{Remark}

\title{Explicit covers and rigid obstructions for simplicial complexity:\\
 from spheres to surfaces}

\author{M. Perez\footnote{Escuela Superior de F\'isica
y Matem\'aticas del Instituto Polit\'ecnico Nacional, Edificio 9, Unidad
Profesional Adolfo L\'opez Mateos, 07300 Mexico City, Mexico.}, and
J. Flores\footnotemark[1]}

\date{}

\begin{document}
\sloppy
\maketitle

\begin{abstract}
The strict simplicial complexity $\SCst(K)$, Gonz\'alez's simplicial
complexity at subdivision level zero, is one less than the least number of contiguity domains covering the
ordered product $K\times K$.  We bound it from both sides.  For upper
bounds, every product of two vertex stars is a domain, which gives
$\SCst(K)\le\rho(K)-1$ for an exactly computable set-cover number
$\rho(K)$; the bound is sharp on every even-dimensional simplex boundary.  A
GPU annealer that evolves whole systems of contiguity chains finds
three-domain covers of $\partial\Delta^3,\dots,\partial\Delta^7$ and of the
cube boundary.  Seeded by star covers, it descends from $24$ to $5$ domains
on the icosahedron boundary, and an exact SAT completion of its stalled
state gives $\SCst\le3$ there.  For lower bounds, a subcomplex of $K\times K$
is rigid if its two projections are joined by no contiguity chain of any
length; $k$ domains can cover $K\times K$ only if $k$ classes of facets
avoid every rigid subcomplex, a finite SAT problem.  Rigid subcomplexes
arise from hollow triangles over essential empty triangles and, on
complexes without dominated vertices, from the graph of every
order-preserving simplicial self-map other than the identity.  The
resulting chain-length-free bounds give
$\SCst(\partial\Delta^{n+1})=2>1=\TC(S^n)$ for $n=3,5,7,9$, $\SCst=2>1$ for
the five-vertex M\"obius band, $\SCst\ge3>2$ for the seven-vertex torus, and
$4\le\SCst\le5$ and $4\le\SCst\le6$ for two triangulations of
$\mathbb{RP}^2$.  They also sort the annealer's stalls: on these two
triangulations of $\mathbb{RP}^2$, the Klein bottle and $\partial\Delta^4$,
almost every uncovered facet closes a rigid subcomplex in every domain,
while the icosahedron stall was a search failure.
\end{abstract}



\section{Introduction}\label{sec:intro}

Farber's topological complexity $\TC(X)$ \cite{F03} is the least number of
continuous motion-planning rules on a space $X$, minus one (we use the
reduced convention throughout, one less than in \cite{F03,FTY03}, so a
cover by $k$ domains certifies $\TC\le k-1$).  Gonz\'alez's simplicial complexity \cite{G18} discretises
it.  On the ordered product $K\times K$ of a triangulation $K$ of $X$,
continuity is replaced by contiguity of simplicial maps
(Section~\ref{sec:back}), and a motion planner on a subcomplex becomes a
chain of contiguous simplicial maps joining the two projections.  The strict
invariant $\SCst(K)$, which is Gonz\'alez's $\SC^0$ \cite{G18} in the
notation of \cite{OLGB21}, counts such planners on $K\times K$ itself,
without subdividing.  It is finite and checkable and satisfies
$\SCst(K)\ge\TC(\|K\|)$; the analogous count on iterated barycentric
subdivisions of $K\times K$ becomes equal to $\TC$ \cite{G18}.  At subdivision level $0$ it is an
invariant of the triangulation, and it is the quantity bounded by the covers
that the randomized greedy algorithm of Ortiz, Lara, Gonz\'alez and Borat
\cite{OLGB21} constructs; values were known for the circle and the wedge of
two circles \cite{OLGB21}.

\paragraph{Results.}  We compute or bracket $\SCst$ at level $0$ for the
simplex boundaries up to $S^9$, four triangulations of $S^2$ (the simplex,
the cube, the icosahedron and the dodecahedron boundaries), the
vertex-minimal M\"obius band, the $7$-vertex torus, an $8$-vertex Klein
bottle and four triangulations of $\mathbb{RP}^2$.  The main statements are:
\begin{itemize}
\item every product $D(v,w)$ of two vertex stars is a contiguity domain, so
$\SCst(K)\le\rho(K)-1$ for an exactly computable set-cover number $\rho(K)$
(Theorem~\ref{thm:prod}, Corollary~\ref{cor:rho});
\item on a complex without dominated vertices the graph $\Gamma_w$ of every
non-decreasing simplicial self-map $w\ne\mathrm{id}$ is rigid, that is, its
two projections are joined by no contiguity chain of any length
(Theorem~\ref{thm:graphs});
\item $\SCst(\partial\Delta^{n+1})=2>1=\TC(S^n)$ for $n=3,5,7,9$;
$\SCst=2>1=\TC$ for the five-vertex M\"obius band; $3\le\SCst\le10$ for the
Cs\'asz\'ar torus ($\TC=2$); $4\le\SCst\le5$ for the hemi-icosahedron and
$4\le\SCst\le6$ for the hemicube ($\TC(\mathbb{RP}^2)=3$); these lower bounds
are computer-assisted (Theorem~\ref{thm:floors});
\item $\SCst(\text{icosahedron boundary})\le3$, so $2\le\SCst\le3$
(Proposition~\ref{prop:icosa});
\item on the hemi-icosahedron, the M\"obius band and the Cs\'asz\'ar torus no
star rectangle $D(v,w)$ can be enlarged by a single facet of $K\times K$ and
remain a domain (Theorem~\ref{thm:maximal}, computer-verified).
\end{itemize}

\paragraph{Method.}  Upper bounds are explicit covers.
Section~\ref{sec:star} gives the star cover, which needs no search; on
$\partial\Delta^{n+1}$ it gives $\SCst\le2$ for all $n$, which is sharp for
even $n$ because $\TC(S^n)=2$.  Section~\ref{sec:annealer} describes a
chain-first annealer, which anneals whole systems of $k$ contiguity chains
on a GPU under an energy that grades every uncovered facet by the defect
count of its best chain.  It finds three-domain covers of $S^2$ through
$S^6$ and of the cube boundary, where $\rho=8$.  Section~\ref{sec:descent}
records descents in $k$ seeded by star covers on the surfaces.

Lower bounds come from rigid subcomplexes (Section~\ref{sec:rigid}).  No
domain contains a rigid subcomplex, so if no $k$ classes of facets avoid a
family of rigid subcomplexes, then $k$ domains cannot cover $K\times K$
(Theorem~\ref{thm:bound}).  The families come from two sources: hollow
triangles whose two projections have different winding, and the graphs
$\Gamma_w$.  The second source contains the rigid spheres that give
$\SCst(\partial\Delta^{n+1})=2$ for odd $n\le9$.  The bounds hold for chains
of every length, which a failed search or a bounded-length solver
infeasibility cannot provide.

The two halves meet in Section~\ref{sec:stall}.  For every stalled
annealer state we test whether some domain can absorb an uncovered facet.
On the hemi-icosahedron, the hemicube, the Klein bottle and
$\partial\Delta^4$ almost no uncovered facet can be absorbed by any domain;
on the icosahedron every uncovered facet can, and an exact SAT decision
(Section~\ref{sec:exact}) completes the stalled state to a $4$-domain cover.

\paragraph{Status of the statements.}  Three kinds of statement appear:
\emph{certificates} (explicit covers, re-checked by a verifier),
\emph{computer-assisted theorems} (unsatisfiable SAT instances built from
proved rigid families, or witness lists re-verified from the definitions),
and \emph{observations} (plateaus and stalls of the search).
Section~\ref{sec:results} collects the values, Section~\ref{sec:repro}
describes the certificates and the auditor, and Section~\ref{sec:concl}
lists open questions.

\section{Background and notation}\label{sec:back}

Throughout, $K$ is a finite, nonempty, connected simplicial complex on a
totally ordered vertex set $V$; a \emph{facet} is a maximal simplex, and
$\operatorname{star}(v)$ denotes the set of facets containing $v$.  The
\emph{ordered product} $K\times K$ is the simplicial complex whose simplices
are the nonempty chains $S\subseteq V\times V$ (for the product order) with
$\pi_1(S)\in K$ and $\pi_2(S)\in K$.  Its realisation is the staircase
triangulation of $\|K\|\times\|K\|$ \cite[Ch.~II]{ES52}, and the projections $\pi_1,\pi_2$ are
simplicial.  A \emph{cell} is a pair $(s,t)$ of facets of $K$; the facets
of $K\times K$ are the maximal chains (monotone lattice paths) in $s\times
t$ over all cells, and when $K$ is pure of dimension $p$ there are
$\binom{2p}{p}$ of them over each cell.  For example, if $K=\partial\Delta^2$
is the circle on $\{0,1,2\}$, then over the cell $(\{0,1\},\{1,2\})$ lie the
two facets $\{(0,1),(0,2),(1,2)\}$ and $\{(0,1),(1,1),(1,2)\}$, and
$K\times K$ has $9\cdot2=18$ facets.  Two simplices of $K\times K$ lie in a
common facet if and only if their union is a simplex.  We write $C_m$
($m\ge3$) for the cycle graph on $y_0,\dots,y_{m-1}$, and a \emph{closed
walk} of length $m$ in a complex is a simplicial map from $C_m$.

Two simplicial maps $\varphi,\varphi':J\to K$ are \emph{contiguous},
$\varphi\sim\varphi'$, if $\varphi(\sigma)\cup\varphi'(\sigma)\in K$ for
every facet $\sigma$ of $J$; a \emph{contiguity chain} with $L$ links is a
sequence $\varphi=\varphi_0\sim\cdots\sim\varphi_L=\varphi'$.  A subcomplex
$D\subseteq K\times K$ is a \emph{contiguity domain} if $\pi_1|_D$ and
$\pi_2|_D$ are joined by a contiguity chain of simplicial maps $D\to K$, and
\[
\SCst(K)=\bigl(\text{least number of domains covering }K\times K\bigr)-1
\;\ge\;\TC(\|K\|)
\]
\cite{G18,OLGB21}; here $\SCst$ is the strict invariant of
\cite[Sec.~1]{OLGB21}.  A cover by $k$ domains is a system of $k$
piecewise-linear motion planners.  All statements here are at subdivision
level $0$, and coverage figures count facets of $K\times K$.

Two elementary facts are used repeatedly.  \emph{Heredity}: a subcomplex of
a domain is a domain (restrict the chain; each facet of the subcomplex lies
in a facet of the domain).  \emph{Facet classes}: if $D_1,\dots,D_k$ cover
$K\times K$ and $\mathcal F_j$ is the set of facets of $K\times K$ lying in
$D_j$, then $\gen{\mathcal F_j}\subseteq D_j$ is a domain and the
$\mathcal F_j$ cover all facets.  So the existence of a $k$-cover is a
question about $k$ sets of facets, which may overlap but can be made
disjoint.

\section{Upper bounds without search: the star cover}\label{sec:star}

Write $c_x$ for the constant map to the vertex $x$.

\begin{theo}[product domains]\label{thm:prod}
For vertices $v,w$ of $K$ let $D(v,w)\subseteq K\times K$ be generated by
the facets over cells $(s,t)$ with $v\in s$, $w\in t$.  Then $D(v,w)$ is a
contiguity domain: if $v=v_0,v_1,\dots,v_m=w$ is an edge path of $K$, then
$\pi_1\sim c_{v_0}\sim c_{v_1}\sim\cdots\sim c_{v_m}\sim\pi_2$ on $D(v,w)$.
\end{theo}

\begin{proof}
A facet $\sigma$ over $(s,t)$ has $\pi_1(\sigma)=s$ and $\pi_2(\sigma)=t$
(a staircase path visits every row and column) and $c_x(\sigma)=\{x\}$, so
all maps are simplicial.  For consecutive maps, $s\cup\{v\}=s$ since $v\in
s$, each $\{v_i,v_{i+1}\}$ is an edge, and $\{w\}\cup t=t$.
\end{proof}

\begin{corollary}\label{cor:rho}
Let $\rho(K)$ be the least size of a set $R\subseteq V\times V$ such that
every cell $(s,t)$ has some $(v,w)\in R$ with $v\in s$, $w\in t$.  Then
$\SCst(K)\le\rho(K)-1$.  For every $n\ge1$, $\rho(\partial\Delta^{n+1})=3$;
hence $\SCst(\partial\Delta^{n+1})\le2$, with equality $=\TC(S^n)=2$ for
even $n\ge2$.
\end{corollary}

\begin{proof}
The first claim follows from Theorem~\ref{thm:prod}, since every facet of
$K\times K$ lies over some cell.  On $\partial\Delta^{n+1}$ with facets
$s_a=V\setminus\{a\}$, the pairs $(0,0),(1,1),(2,2)$ cover: $(s_a,s_b)$ is
missed by $(v,v)$ iff $v\in\{a,b\}$, and $\{0,1,2\}\not\subseteq\{a,b\}$.
Two pairs $(v,w),(v',w')$ never suffice, because both miss $(s_v,s_{w'})$.
Equality for even $n$ follows from $\TC(S^n)=2$ \cite{F03} and
$\SCst\ge\TC$.
\end{proof}

\begin{lemma}[counting bounds]\label{lem:rows}
Let $\tau$ be the least size of a set of vertices meeting every facet,
$\delta=\max_v|\operatorname{star}(v)|$, and $F$ the number of facets.  Then
$\lceil\tau F/\delta\rceil\le\rho(K)\le\tau^2$.
\end{lemma}

\begin{proof}
A transversal $S$ of size $\tau$ gives the cover $S\times S$.  For the lower
bound fix an optimal $R$ and a facet $s$.  The set $W_s=\{w:\exists v\in
s,\ (v,w)\in R\}$ meets every facet $t$, because the cell $(s,t)$ is
covered, so $|W_s|\ge\tau$.  Counting incidences,
$|R|\,\delta\ge\sum_{(v,w)\in R}|\operatorname{star}(v)|\ge\sum_s|W_s|\ge\tau F$.
\end{proof}

The lower bound of the lemma is attained on the simplex boundaries
(Corollary~\ref{cor:rho}), the wedge and the M\"obius band ($\rho=4$), the
hemi-icosahedron ($\rho=6$) and the icosahedron ($\tau=6$, $\delta=5$,
$F=20$, so $\rho\ge24$), so on these complexes an explicit rectangle set of
that size determines $\rho$.
Table~\ref{tab:rho} lists $\rho$ for the complexes of this paper; every
value is the exact optimum of the set-cover instance (CP-SAT \cite{ORTools},
status OPTIMAL, at most $12$\,s).  Exactness matters: greedy set cover
needs $27$ rectangles on the icosahedron against $24$.

\begin{table}[H]
\centering\footnotesize
\setlength{\tabcolsep}{4pt}
\begin{tabular}{@{}lrrrrrrr@{}}
\toprule
complex $K$ & $|V|$ & $F$ & $\tau$ & $\delta$ &
$\lceil\tau F/\delta\rceil$ & $\rho(K)$ & $\tau^2$ \\
\midrule
$\partial\Delta^{n+1}$, $n\ge1$ & $n{+}2$ & $n{+}2$ & $2$ & $n{+}1$ & $3$ & $3$ & $4$\\
wedge $S^1\vee S^1$ & $4$ & $5$ & $2$ & $3$ & $4$ & $4$ & $4$\\
M\"obius band $K_M$ & $5$ & $5$ & $2$ & $3$ & $4$ & $4$ & $4$\\
cube boundary & $8$ & $12$ & $3$ & $6$ & $6$ & $8$ & $9$\\
hemi-icosahedron ($\mathbb{RP}^2$) & $6$ & $10$ & $3$ & $5$ & $6$ & $6$ & $9$\\
hemicube ($\mathbb{RP}^2$) & $7$ & $12$ & $3$ & $6$ & $6$ & $7$ & $9$\\
hemi-octahedron ($\mathbb{RP}^2$) & $9$ & $16$ & $3$ & $6$ & $8$ & $9$ & $9$\\
hemi-dodecahedron ($\mathbb{RP}^2$) & $10$ & $18$ & $4$ & $9$ & $8$ & $15$ & $16$\\
Cs\'asz\'ar torus & $7$ & $14$ & $4$ & $6$ & $10$ & $11$ & $16$\\
Klein bottle (Cervone \cite{Ce94}) & $8$ & $16$ & $4$ & $7$ & $10$ & $12$ & $16$\\
icosahedron boundary & $12$ & $20$ & $6$ & $5$ & $24$ & $24$ & $36$\\
dodecahedron boundary & $20$ & $36$ & $6$ & $8$ & $27$ & $36$ & $36$\\
\bottomrule
\end{tabular}
\caption{The rectangle number $\rho(K)$ and the bounds of
Lemma~\ref{lem:rows}.  Facet lists are in \texttt{src/gpu\_sc.py} of the
source repository; optimal rectangle sets are in
\texttt{results/*\_starcover*/summary.json}.  The square faces of the cube
are split by one diagonal each and the pentagons of the dodecahedron are
fanned from a vertex.  The hemicube and hemi-octahedron are subdivided to be
simplicial (a centre vertex per square, a midpoint per edge); the pentagons
of the hemi-dodecahedron are fanned.  The four triangulations of $S^2$ have
$\rho=3,8,24,36$, so $\rho$ depends on the triangulation, not only on the
space.}
\label{tab:rho}
\end{table}

\paragraph{The star cover as a method} (\texttt{src/star\_cover.py},
\texttt{src/rect\_cover.py}).  Solve the set cover exactly, chain each
rectangle along a shortest edge path, assign every facet to the first
rectangle containing it, shorten each chain by the greedy reduction of
\cite[Alg.~2]{OLGB21}, and verify.  The method cannot fail and takes
seconds (the icosahedron in $1.3$\,s with greedy set cover and $12$\,s with
the exact one; $\partial\Delta^9$, with $1\,287\,000$ facets, by three
domains).  It returns a complete verified cover, the bound $\rho-1$, and
chains that seed the annealer.  The gap between $\rho-1$ and $\SCst$ comes
from restricting to product-shaped domains with constant chains.  It is
large on sparse complexes (cube: $\rho-1=7$ against $\SCst=2$; icosahedron:
$23$ against at most $3$) and vanishes on every simplex boundary
$\partial\Delta^{n+1}$ with $2\le n\le9$ (Corollary~\ref{cor:rho} and
Theorem~\ref{thm:floors}).

\section{Upper bounds by search: the chain-first annealer}\label{sec:annealer}

\paragraph{Growing domains one facet at a time.}  The pipeline of
\cite{OLGB21} grows one domain at a time by budgeted random walks
(\textsc{LocalSearch}).  Our GPU port of it on
$\partial\Delta^3\times\partial\Delta^3$ ($96$ facets, walk budget $40000$,
$16$ replicas) froze every replica at a \emph{near-rectangle} of $30$ to
$54$ facets, and $65$ consecutive retry passes over a frozen front, which
took $82\%$ of the wall time, produced no extension.  An earlier C++/MPI
version of ours on $96$ cores never completed a cover in $120$-hour jobs.
These are our runs, not measurements reported in \cite{OLGB21}.  Two
features of the method are consistent with this record: a failed extension
returns no graded distance to success, while contiguity is a hard
constraint; and growth is sequential and irrevocable, so it can stop at
fronts where no tested extension succeeds.  Theorem~\ref{thm:maximal}
below shows that on three of the surfaces a star rectangle cannot be
extended by any single facet.

\paragraph{The annealer.}  A \emph{state} is a $k$-tuple
$C=(C_1,\dots,C_k)$ of sequences $C_j=(\varphi^j_0,\dots,\varphi^j_{c_{\max}})$
of vertex assignments $V\times V\to V$ with $\varphi^j_0=\pi_1$ and
$\varphi^j_{c_{\max}}=\pi_2$; the interior entries are arbitrary.  A facet
$\sigma$ is \emph{covered by} $C_j$ when every link
$\varphi^j_i\sim\varphi^j_{i+1}$ is contiguous on $\sigma$.  On the
subcomplex $\mathrm{Dom}(C_j)$ generated by the facets covered by $C_j$,
the entries are simplicial and consecutive entries are contiguous (each
image of a covered facet lies in the union of the images of an adjacent
link, which is a simplex of $K$), so $C_j$ restricted to $\mathrm{Dom}(C_j)$
is a contiguity chain and $\mathrm{Dom}(C_j)$ is a domain.  Let
$b_j(\sigma)$ be the number of links of $C_j$ that are not contiguous on
$\sigma$, $m(\sigma)=\min_j b_j(\sigma)$, $\mathrm{pen}(0)=0$ and
$\mathrm{pen}(b)=1+\beta b$ for $b>0$.  The energy is
\[
E=\sum_\sigma \omega_\sigma\,\mathrm{pen}(m(\sigma))+\gamma\sum_j\sum_\sigma b_j(\sigma),
\]
with facet weights $\omega_\sigma=1$ unless a guided rung raises them.  The
first sum vanishes exactly at a cover and grades every uncovered facet by
its best chain, which is the graded signal the greedy method lacks.  The
second sum penalises defects in all chains and keeps chains from drifting
($\gamma=0$ made long chains fail in the sphere runs).  \emph{Coverage},
the reported quantity, is the number of facets with $m(\sigma)=0$.  The
stopping criterion is full coverage followed by verification, never $E=0$.

A population of $S$ states anneals in parallel on the GPU under a
Metropolis rule \cite{KGV83}.  States are initialised by random sweeps that
flip $\pi_1$ towards $\pi_2$.  A proposal changes the value of one vertex
in one chain, at a single interior position or over a block of two to
eight consecutive positions (\emph{track} moves, probability
$p_{\mathrm{track}}$), copying from a neighbouring map with probability
$1/2$ and otherwise drawing a fresh vertex.  The surface campaigns add
compound \emph{stretch}, \emph{collapse} and \emph{seam} proposals that
rewrite several columns of a block at once.  The temperature decreases
geometrically within each epoch; every second epoch the weakest quarter of
the population is overwritten by the best states and a few slots are
re-initialised.  The proposal kernel is not symmetric, and we claim no
convergence guarantee.  When the incumbent covers every facet, each facet
is assigned to the first chain covering it, empty domains are discarded,
chains are reduced, and the system is verified map by map, link by link and
facet by facet; the output is a certificate.  Algorithm~\ref{alg:anneal}
summarises the loop.

\begin{algorithm}[H]
\caption{\textsc{ChainAnneal}($K$, $k$, $c_{\max}$; $S$, $\beta$, $\gamma$,
$p_{\mathrm{track}}$).}
\label{alg:anneal}
\begin{algorithmic}[1]
\State initialise $S$ states by random sweeps from $\pi_1$ to $\pi_2$;
$\mathcal C^*\gets$ incumbent (highest coverage)
\For{epoch $=1,2,\dots$}
\ForAll{states, in parallel on the GPU}
\State $M$ Metropolis steps at temperature $T_0\to T_1$: pick a chain $j$
and a vertex $x$, propose a new value at one position or over a block of
$2$--$8$ positions (copied from a neighbour w.p.\ $1/2$, else fresh), accept
w.p.\ $\min\{1,e^{-\Delta E/T}\}$
\EndFor
\State update $\mathcal C^*$ by unweighted coverage
\If{$\mathcal C^*$ covers every facet}
\State partition facets to chains, discard empty domains, reduce each chain
\cite[Alg.~2]{OLGB21}, verify; \Return the certificate ($\SCst\le k-1$)
\EndIf
\State every second epoch: overwrite the weakest quarter with the best
states, re-initialise a few slots
\EndFor
\end{algorithmic}
\end{algorithm}

\paragraph{Ladders.}  Runs are organised in \emph{rungs}, successive runs
each started from the result of the previous one.  A \emph{warm} rung
initialises half the population from the previous incumbent (chains reduced
on their domains and padded with $\pi_2$); a \emph{guided} rung also raises
the weight of each facet the incumbent misses by $g$ after each epoch
(capped at $g_{\max}$, decaying towards $1$ by a factor $0.7$).  A rung has
a fixed overhead $t_0\approx40$\,s, and its gain is front-loaded (an
empirical fit on the dodecahedron gives $\Delta\mathrm{cov}(t)\approx
c\,t^{0.8}$); maximising $c\,t^{0.8}/(t+t_0)$ gives rungs of length
$4t_0\approx2.7$\,min.  The surface campaigns also use a
\emph{$k$-descent} ladder: a complete $k$-cover seeds the $(k{-}1)$-search
with its $k-1$ largest domains (Section~\ref{sec:descent}).

\begin{table}[H]
\centering\footnotesize
\setlength{\tabcolsep}{3.5pt}
\begin{tabular}{@{}llr@{\hspace{8pt}}crl>{\raggedright\arraybackslash}p{3.1cm}@{}}
\toprule
target & $K$ & facets & $k$ & coverage & wall time & certified domains \\
\midrule
$S^1$ & $\partial\Delta^2$ & $18$ & $2$ & $18/18$ & $7.9$\,s & $9+9$; $\SCst=1=\TC$\\
$S^2$ & $\partial\Delta^3$ & $96$ & $3$ & $96/96$ & $22$\,s & $54+29+13$; $\SCst=2=\TC$\\
$S^3$ & $\partial\Delta^4$ & $500$ & $3$ & $500/500$ & $153$\,s & $309+130+61$; $\SCst\le2$\\
$S^4$ & $\partial\Delta^5$ & $2520$ & $3$ & $2520/2520$ & $13$\,min & $1634+653+233$; $\SCst=2=\TC$\\
$S^5$ & $\partial\Delta^6$ & $12348$ & $3$ & $12348/12348$ & $34$\,min & $7061+4239+1048$; $\SCst\le2$\\
$S^6$ & $\partial\Delta^7$ & $59136$ & $3$ & $59136/59136$ & ${\approx}2.5$\,h$^\S$ & $44909+11451+2776$; $\SCst=2=\TC$\\
\midrule
$S^1\vee S^1$ & wedge & $50$ & $3$ & $50/50$ & $41$\,s & $18+16+16$; $\SCst=2$\\
$S^2$ & cube boundary & $864$ & $3$ & $864/864$ & ${\approx}24$\,min$^\dagger$ & $482+236+146$; $\SCst=2=\TC$\\
M\"obius & $K_M$ & $150$ & $4$ & $150/150$ & $17$\,s$^\ddagger$ & $54+37+34+25$; $\SCst\le3$\\
\midrule
$S^2$ & $\partial\Delta^3$ & $96$ & $2$ & $84/96$ & plateau & two-rectangle value\\
$S^3$ & $\partial\Delta^4$ & $500$ & $2$ & $460/500$ & plateau & two-rectangle value ($8$ runs)\\
$S^4$ & $\partial\Delta^5$ & $2520$ & $2$ & $2380/2520$ & $12$\,min & two-rectangle value\\
\bottomrule
\end{tabular}
\caption{Chain-first annealer on one Tesla T4 unless noted; $k$ chains
export at most $k$ nonempty domains.  Parameters: $c_{\max}=20$--$40$,
$S=4096$ (circle $1024$), $5000$--$20000$ moves per epoch, $\beta=5$,
$\gamma=0.5$, $p_{\mathrm{track}}=0.3$ (M\"obius: $\beta=1$,
$p_{\mathrm{track}}=0.5$).  $^\dagger$Three rungs: cold $303$\,s
($662/864$), warm $814$\,s ($820$), guided $304$\,s (complete); a guided
cold rung alone stalled at $858$.  $^\ddagger$Warm from a $k=5$ cover
($36$\,s); eight $k=3$ runs stalled at $148/150$ (Section~\ref{sec:stall}).
$^\S$Eight-rung ladder on an RTX PRO 6000 ($S=64$ cold, $S=256$ guided).  The
$k=2$ plateaus equal the number of facets covered by the best pair of star
rectangles ($14$, $23$ and $34$ cells of $6$, $20$ and $70$ facets); the
domains realising them are often not rectangles.}
\label{tab:annealer}
\end{table}

\paragraph{The dodecahedron.}  On the dodecahedron boundary ($20$
vertices, $36$ triangles, $7776$ product facets; $3.2\%$ of vertex triples
are simplices, and the dual graph has diameter $7$) a $20$-hour campaign at
$k=3$ (five warm escalations of up to an hour, then $172$ rungs of
$2.5$\,min over $12$ GPU-hours) rose from $644$ to $2602$ facets ($33.5\%$),
with mean gain per rung falling from $11.7$ to below $1$.  No cover with
fewer than $\rho=36$ domains is known, so the star bound $\SCst\le35$
stands.

\section{Star-seeded descents on the surfaces}\label{sec:descent}

The surface campaigns ran \texttt{annealer\_stretch.py} ($\beta=5$,
$\gamma=0.1$, $c_{\max}=40$, $S=4096$, $5000$ steps per epoch, seed $8233$
unless stated) with the $k$-descent ladder.  Table~\ref{tab:descent}
records first covers, descents and first stalls.  Finding a first cover from
a cold start did not scale: it took $2$\,h on the Klein bottle, and cold
runs on the hemi-octahedron and hemi-dodecahedron ended without a cover.
The star cover replaces this stage.  The hemi-dodecahedron, star-seeded at
$k=14$, was covered in $6689$\,s, where cold runs at $k=15,20,30$ had found
no cover after $3521$, $14673$ and $2863$\,s.  The icosahedron descended
from the $24$-domain star cover to a verified $5$-domain cover in $17.9$\,h
of rungs, every rung from $k=23$ to $k=5$ closing (the $k=5$ rung after a
warm continuation of a $2395/2400$ stall).  Each complex was run with a
single seed and a single lineage, so these are observations, not controlled
comparisons.

\begin{table}[H]
\centering\footnotesize
\setlength{\tabcolsep}{3pt}
\begin{tabular}{@{}l r>{\raggedright\arraybackslash}p{2.9cm}>{\raggedright\arraybackslash}p{2.5cm}>{\raggedright\arraybackslash}p{1.9cm}>{\raggedright\arraybackslash}p{2.4cm}@{}}
\toprule
complex ($|K{\times}K|$) & $\rho$ & first cover & descent (wall) & lowest verified & first stall \\
\midrule
hemi-icosahedron ($600$) & $6$ & cold $k{=}11$: $679$\,s & $k{=}11\to6$ ($0.8$\,h) & $6$ domains & $k{=}5$: $540$--$544/600$\\
hemicube ($864$) & $7$ & cold $k{=}20$: $2470$\,s & $k{=}20\to7$ ($2.0$\,h) & $7$ & $k{=}6$: $836/864$\\
Klein bottle ($1536$) & $12$ & cold $k{=}15$: $7075$\,s & $k{=}15\to10$ ($4.3$\,h) & $10$ & $k{=}9$: $1494/1536$\\
hemi-octahedron ($1536$) & $9$ & star-seeded $k{=}8$: $6298$\,s & --- & $8$ & $k{=}7$: $1432/1536$ (log)\\
hemi-dodecahedron ($1944$) & $15$ & star-seeded $k{=}14$: $6689$\,s & $k{=}14\to7$ ($12.3$\,h) & $7$ & $k{=}6$: $1922/1944$\\
icosahedron ($2400$) & $24$ & star-seeded $k{=}23$: $1545$\,s & $k{=}23\to5$ ($17.9$\,h) & $5$; then $4$ (Section~\ref{sec:exact}) & $k{=}4$: $2376/2400$ ($3$ runs)\\
dodecahedron ($7776$) & $36$ & none below $\rho$ & --- & star only & $k{=}3$: $2602/7776$\\
Cs\'asz\'ar torus ($1176$) & $11$ & none below $\rho$ & --- & star only & $k{=}11$ cold: $563/1176$\\
\bottomrule
\end{tabular}
\caption{First covers, descents and first stalls.  Sources:
\texttt{results/rp2\_k11}\dots\texttt{rp2\_k6},
\texttt{hemicube\_k20}\dots\texttt{k7}, \texttt{klein\_k15\_8233\_v2}\dots\texttt{klein\_k10},
\texttt{hemi\_octa\_k8}, \texttt{hemi\_dodeca\_k14}\dots\texttt{k7},
\texttt{icosa\_k23}\dots\texttt{icosa\_k5\_v2}; stalls in \texttt{rp2\_k5\_*},
\texttt{hemicube\_k6\_v2}, \texttt{klein\_k9\_v4}, \texttt{hemi\_dodeca\_k6\_v2},
\texttt{icosa\_k4\_v3}, \texttt{dodeca\_*}.  A stall is the best coverage of
the recorded runs at that $k$; Section~\ref{sec:stall} decides which stalls
are obstructions.}
\label{tab:descent}
\end{table}

\section{Lower bounds for chains of every length: rigid subcomplexes}\label{sec:rigid}

A failed search, or even a solver infeasibility for chains of at most $L$
links, says nothing about longer chains (a shortest chain has distinct
maps, so its length is bounded a priori, but only by the number of vertex
maps, far beyond reach).  The bounds of this section hold
for chains of every length.  The notion of a rigid subcomplex and the
winding argument of Section~\ref{sec:src1} were first written up in the
authors' note \cite{PM26b}; the arguments needed here are given in full.

\begin{definition}\label{def:rigid}
A subcomplex $Z\subseteq K\times K$ is \emph{rigid} if $\pi_1|_Z$ and
$\pi_2|_Z$ are joined by no contiguity chain of simplicial maps $Z\to K$, of
any length; equivalently, $Z$ is not a contiguity domain.  For a family
$\mathcal Z$ of subcomplexes, a \emph{$\mathcal Z$-avoiding $k$-cover} is a
family $\mathcal F_1,\dots,\mathcal F_k$ of sets of facets of $K\times K$
covering all facets such that no $\gen{\mathcal F_j}$ contains a member of
$\mathcal Z$; $\chi(\mathcal Z)$ is the least such $k$.
\end{definition}

\begin{theo}[rigid-subcomplex bound]\label{thm:bound}
If $\mathcal Z$ is a family of rigid subcomplexes, then
$\SCst(K)+1\ge\chi(\mathcal Z)$.  Deciding whether a $\mathcal Z$-avoiding
$k$-cover exists is a finite SAT instance (variables ``facet $f$ is in
class $j$'' and ``simplex $z$ is absent from class $j$'').
\end{theo}

\begin{proof}
The facet classes of a cover by $\SCst(K)+1$ domains
(Section~\ref{sec:back}) generate domains, and by heredity a domain
contains no rigid subcomplex.
\end{proof}

The bound is tight in principle: $K\times K$ has $|V|^2$ vertices, so if
$\mathcal Z$ is the family of all rigid subcomplexes, then a
$\mathcal Z$-avoiding class $\gen{\mathcal F_j}$ is itself not rigid, hence
a domain, and $\chi(\mathcal Z)=\SCst(K)+1$.  In practice one uses small
families that can be enumerated.

\subsection{Source 1: closed walks with different contiguity classes}\label{sec:src1}

Let $u$ be a $\Zt$-valued $1$-cocycle of $K$ (with $u=0$ on degenerate
edges $\{a,a\}$), and for a closed walk $f:C_m\to K$ let
$\nu_u(f)=\sum_i u(f(y_i)f(y_{i+1}))$ (indices mod $m$).

\begin{lemma}[contiguity preserves winding]\label{lem:winding}
If $f,f':C_m\to K$ are contiguous, then $\nu_u(f)=\nu_u(f')$.
\end{lemma}

\begin{proof}
Put $a=f(y_i)$, $b=f(y_{i+1})$, $a'=f'(y_i)$, $b'=f'(y_{i+1})$.
Contiguity on the edge $y_iy_{i+1}$ makes $\{a,b,a',b'\}$ a simplex, so
$\{a,b,b'\}$ and $\{a,a',b'\}$ are simplices and the cocycle relation gives
$u(ab)+u(bb')=u(ab')=u(aa')+u(a'b')$.  Summing over $i$, the terms
$u(f(y_i)f'(y_i))$ appear twice and cancel, leaving $\nu_u(f)=\nu_u(f')$.
\end{proof}

A \emph{hollow triangle} $Z\subseteq K\times K$ consists of three pairwise
comparable grid points whose pairs are edges of $K\times K$ and whose
triple is not a simplex.  If a contiguity chain joined $\pi_1|_Z$ and
$\pi_2|_Z$, then $\pi_1\circ Z$ and $\pi_2\circ Z$ would be contiguous
closed walks; so by Lemma~\ref{lem:winding}, $Z$ is rigid as soon as
$\nu_u(\pi_1\circ Z)\ne\nu_u(\pi_2\circ Z)$.  The basic case is the
\emph{slice} $T\times\{c\}$ over an \emph{essential empty triangle} $T$ of
$K$ (three edges of $K$ spanning no face, with $u(\partial T)=1$), on which
$\pi_1$ winds once and $\pi_2$ is constant.

More generally, let $Z$ be the subcomplex of $K\times K$ traced by a
closed walk $g:C_m\to K\times K$.  A chain on $Z$ composed with $g$ is a
chain of maps $C_m\to K$ (precomposition with a simplicial map preserves
contiguity), so $Z$ is rigid whenever $\pi_1\circ g$ and $\pi_2\circ g$
lie in different \emph{contiguity classes} of maps $C_m\to K$.  We computed
these classes exactly (union--find over all $|V|^m$ vertex maps).  They are
finer than homotopy classes.  The $3$-walks of the Klein bottle fall into
$11$ classes and its $4$-walks into $6$; the torus separates the two
integral windings ($7$ classes of $3$-walks); the hemi-octahedron has $25$
classes of $3$-walks and $2$ of $4$-walks.  The hemi-icosahedron has
exactly $2$ classes (the $\Zt$ winding) for every length from $3$ to $6$,
and the icosahedron and cube boundaries have a single class for every
length up to $5$ and $6$ respectively, so no short closed walk is rigid
there.

\subsection{Source 2: rigid graphs over minimal complexes}

A vertex $v$ is \emph{dominated} if some $v'\ne v$ lies in every facet
containing $v$; so $K$ has no dominated vertex iff
$\bigcap\{\sigma\text{ facet}:v\in\sigma\}=\{v\}$ for all $v$.  These are
the \emph{minimal} complexes of Barmak and Minian \cite{BM12}.  Every
closed surface triangulation is minimal (the facets at $v$ are the
triangles around the link cycle of $v$, which has length at least $3$), and
so are $\partial\Delta^{n+1}$ ($n\ge1$) and the five-vertex M\"obius band.
More generally, a complex without free faces is minimal: if $v'$ dominates
$v$ and $\eta$ is a maximal simplex of the link of $v$ not containing $v'$,
then $\eta\cup\{v\}$ is a free face of $\eta\cup\{v,v'\}$.

The following lemma is a special case of \cite[Prop.~2.7]{BM12}, which
shows that on a minimal complex the identity is the only map in its
contiguity class; we include the one-line proof.

\begin{lemma}\label{lem:idalone}
If $K$ has no dominated vertex and $\varphi:K\to K$ is simplicial and
contiguous to $\mathrm{id}_K$, then $\varphi=\mathrm{id}_K$.
\end{lemma}

\begin{proof}
For a facet $\sigma$, $\sigma\cup\varphi(\sigma)\in K$ and $\sigma$ is
maximal, so $\varphi(\sigma)\subseteq\sigma$; hence
$\varphi(v)\in\bigcap_{\sigma\ni v}\sigma=\{v\}$.
\end{proof}

\begin{theo}[rigid graphs]\label{thm:graphs}
Let $K$ have no dominated vertex and let $w:V\to V$ be non-decreasing for
the vertex order and simplicial ($w(\sigma)\in K$ for all $\sigma\in K$),
with $w\ne\mathrm{id}$.  Then
\[
\Gamma_w=\bigl\{\{(v,w(v)):v\in\sigma\}:\sigma\in K\bigr\}
\]
is a subcomplex of $K\times K$, $\pi_1|_{\Gamma_w}$ is an isomorphism onto
$K$ under which $\pi_2|_{\Gamma_w}$ becomes $w$, and $\Gamma_w$ is rigid.
The same holds for the transpose $\Gamma_w^{\mathsf T}=\{\{(w(v),v):v\in\sigma\}:\sigma\in K\}$.
\end{theo}

\begin{proof}
For $\sigma\in K$, the set $S_\sigma=\{(v,w(v)):v\in\sigma\}$ is a chain
because $w$ is non-decreasing, and $\pi_1(S_\sigma)=\sigma\in K$,
$\pi_2(S_\sigma)=w(\sigma)\in K$ because $w$ is simplicial; so $S_\sigma$ is
a simplex of $K\times K$, and its faces are the $S_{\sigma'}$ with
$\sigma'\subseteq\sigma$.  The bijection $(v,w(v))\mapsto v$ carries
$S_\sigma$ to $\sigma$, hence is an isomorphism $\Gamma_w\to K$ under which
$\pi_1$ becomes $\mathrm{id}_K$ and $\pi_2$ becomes $w$.  A contiguity chain
$\pi_1|_{\Gamma_w}=\psi_0\sim\cdots\sim\psi_L=\pi_2|_{\Gamma_w}$ therefore
becomes a chain $\mathrm{id}=\varphi_0\sim\varphi_1\sim\cdots\sim\varphi_L=w$
of simplicial maps $K\to K$.  By Lemma~\ref{lem:idalone} and induction on
$i$, every $\varphi_i=\mathrm{id}$, so $w=\mathrm{id}$, a contradiction.  For
the transpose, read the chain backwards.
\end{proof}

The hypothesis on dominated vertices cannot be dropped: on the path
$0$--$1$--$2$ the constant map $w\equiv1$ is non-decreasing and simplicial,
and $\pi_1\sim c_1=\pi_2$ on $\Gamma_w=K\times\{1\}$.

\begin{remark}\label{rem:graphs}
(a) For $K=\partial\Delta^{n+1}$ every non-decreasing $w\ne\mathrm{id}$ is
simplicial (it is not onto, so every simplex goes to a proper subset of
$V$), and $\Gamma_w$ is a hollow $n$-sphere in $K\times K$.  These are the
rigid spheres of \cite{PM26b}; Theorem~\ref{thm:graphs} extends them to
every minimal complex.  (b) For $K$ without dominated vertices, the
constants $w\equiv c$ give $\Gamma_c=K\times\{c\}$, so no domain contains a
full row or column of $K\times K$.  (c) A contractible complex without free
faces, such as a triangulation of the dunce hat, has $K\times\{c\}$ rigid
although $\pi_1$ and $\pi_2$ are homotopic on it, so $\SCst(K)\ge1>0=\TC$;
this rigidity is combinatorial, not homotopical.  On the complexes we
tested, however, every $\Gamma_w$ has even mod-$2$ degree or acts
non-trivially on $H^1(K;\Zt)$, so $w$ is not homotopic to the identity
($3891$ maps on the icosahedron, $597$ on the Klein bottle, $181$ on the
hemi-icosahedron, no exception).  (d) We verified by brute force that the
identity is alone in its contiguity class on $\partial\Delta^3$,
$\partial\Delta^4$, the cube, the icosahedron, the hemi-icosahedron, the
hemicube, the hemi-octahedron, the hemi-dodecahedron, the Klein bottle, the
torus and the M\"obius band, whereas on the $2$-simplex all $27$ vertex
maps are contiguous to the identity and on the path $0$--$1$--$2$ four are.
\end{remark}

\subsection{The bounds}

For each complex we built the family $\mathcal Z$ of all rigid closed walks
of length at most $4$ (contiguity classes computed exactly) together with
all $\Gamma_w$ and their transposes, and decided the avoidance instances
with CaDiCaL \cite{CaDiCaL} through PySAT \cite{PySAT}
(\texttt{scripts/rigid\_analysis/rigid.py}).
Table~\ref{tab:floors} gives the results.  An UNSAT entry at $k$ proves
that $k$ domains do not suffice, for chains of every length.

\begin{table}[H]
\centering\footnotesize
\setlength{\tabcolsep}{3pt}
\begin{tabular}{@{}lrrrrcc@{}}
\toprule
$K$ & rigid walks & rigid graphs & $\TC{+}1$ & lower bound & UNSAT at $k$ & SAT at $k$ \\
\midrule
M\"obius band $K_M$ & $4\,684$ & $150$ & $2$ & $\mathbf 3$ & $2$ & $3$ \\
$\partial\Delta^4$ ($S^3$) & $0$ & $250$ & $2$ & $\mathbf 3$ & $2$ & $3$ \\
hemi-icosahedron & $19\,214$ & $362$ & $4$ & $\mathbf 5$ & $4$ ($28$\,s) & $5$ ($540$\,s) \\
hemicube & $29\,216$ & $900$ & $4$ & $\mathbf 5$ & $4$ ($486$\,s) & $5$ ($411$\,s) \\
Klein bottle & $94\,540$ & $1\,194$ & $5$ & $5$ & $4$ ($911$\,s) & $5$: undecided \\
Cs\'asz\'ar torus & $87\,204$ & $686$ & $3$ & $\mathbf 4$ & $3$ & $4$ \\
hemi-octahedron & $30\,922$ & $2\,056$ & $4$ & $4$ & $3$ & $4$ ($1476$\,s) \\
hemi-dodecahedron & $37\,966$ & $2\,042$ & $4$ & $4$ & $3$ & $4$ \\
cube, icosahedron & $0$ & $844$, $7\,782$ & $3$ & $2$ & --- & $2$ \\
\bottomrule
\end{tabular}
\caption{Lower bounds on the number of domains from
Theorem~\ref{thm:bound}, with the families of Sources 1 and 2 (distinct edge
sets of rigid walks of length $\le4$; rigid graphs with transposes).  ``SAT
at $k$'' means that the family admits an avoiding $k$-cover, so it gives
nothing beyond $k-1$; it does not mean that $k$ domains suffice.  Bold
entries exceed $\TC+1$.  On the cube and the icosahedron the family gives
less than $\TC+1$.  The solver runs have no checked proof logs.  The best
known covers are in Table~\ref{tab:surfaces}.}
\label{tab:floors}
\end{table}

For the odd spheres the same test was run on the rigid-sphere family of
Remark~\ref{rem:graphs}(a) at $k=2$ (\texttt{scripts/relaxation/spheres.py}
and \texttt{spheres\_big.py} in the source repository): all $250$ graphs and
transposes for $\partial\Delta^4$ and all $3430$ for $\partial\Delta^6$ (six
SAT solvers agreeing), the $5298$ graphs of maps with at most three values
for $\partial\Delta^8$ (three solvers) and the $15972$ such graphs for
$\partial\Delta^{10}$ (two solvers).  All four instances are
unsatisfiable.  On the circle $\partial\Delta^2$ the instance is
satisfiable, as it must be.

\begin{theo}[computer-assisted]\label{thm:floors}
At subdivision level $0$:
$\SCst(\partial\Delta^{n+1})=2>1=\TC(S^n)$ for $n=3,5,7,9$;
$\SCst(K_M)=2>1=\TC$ for the five-vertex M\"obius band;
$3\le\SCst(\text{Cs\'asz\'ar torus})\le10$, with $\TC=2$ \cite{F03};
$4\le\SCst(\text{hemi-icosahedron})\le5$ and $4\le\SCst(\text{hemicube})\le6$,
with $\TC(\mathbb{RP}^2)=3$ \cite{FTY03};
$4\le\SCst(\text{Klein bottle})\le9$, with $\TC=4$ \cite{CV17}.
\end{theo}

\begin{proof}
Lower bounds: Theorem~\ref{thm:bound} applied to the UNSAT instances above
(Table~\ref{tab:floors} and the rigid-sphere instances), the families being
rigid by the contiguity-class argument of Section~\ref{sec:src1} and by
Theorem~\ref{thm:graphs}; for the Klein bottle the lower bound is $\TC$.
Upper bounds: Corollary~\ref{cor:rho} for the spheres, the star cover for
the torus, and the certificates of Tables~\ref{tab:annealer},
\ref{tab:descent} and~\ref{tab:surfaces} for the rest.  The hemicube upper
bound rests on a certificate checked by the search's own verifier only
(Section~\ref{sec:repro}).
\end{proof}

For $\partial\Delta^4$ this replaces the CP-SAT result of
Section~\ref{sec:sat} (no two-domain cover with chains of at most $48$
links) by a statement for chains of every length.  The M\"obius value had
been established for chains of at most $8$ links \cite{PM26}; it now holds
for chains of every length (a proof by hand is in \cite{PM26b}).  For the
hemi-icosahedron the relaxation that forbids, in every class, all closed
walks of every length with mismatched $\Zt$ winding
(\texttt{scripts/relaxation/h1cegar.py}) is satisfiable at $k=5$, so the
winding argument alone cannot exclude a five-domain cover; whether one
exists is open (Section~\ref{sec:exact}).

\section{Anatomy of a stall}\label{sec:stall}

The families of Section~\ref{sec:rigid} decide whether a stalled state can
be extended one facet at a time.  Let a state have domains
$D_1,\dots,D_k$ and set of uncovered facets $U$.  An uncovered facet $f$ is
an \emph{open door} of $D_j$ if the subcomplex generated by $D_j$ and $f$
contains no member of the family, and $f$ is \emph{dead} if it is an open
door of no $D_j$.  If $f$ is dead, no chain of any length covers $D_j$ and
$f$ together, for any $j$; a search can cover $f$ only after some domain
gives up facets.

\begin{theo}[maximal rectangles; computer-verified]\label{thm:maximal}
Let $K$ be the hemi-icosahedron, the five-vertex M\"obius band or the
Cs\'asz\'ar torus.  For every vertex pair $(v,w)$ and every facet
$f\notin D(v,w)$ of $K\times K$, the subcomplex generated by $D(v,w)$ and $f$
contains a hollow triangle whose two projected $3$-walks lie in different
contiguity classes of maps $C_3\to K$.  Hence no subcomplex generated by
$D(v,w)$ and further facets of $K\times K$ is a domain: every star
rectangle is maximal among domains generated by facets.
\end{theo}

\begin{proof}
There are $36\cdot450$, $25\cdot96$ and $49\cdot960$ pairs $((v,w),f)$, and
we list an explicit witness for each
(\texttt{scripts/rigid\_analysis/}\allowbreak\texttt{certificates/}\allowbreak\texttt{maximal\_rectangles\_*.json}).
Each witness was re-verified by a separate program: every edge of the
triangle lies in a facet of $D(v,w)$ or in $f$, and the two walks lie in
different classes, the classes ($2$, $3$ and $7$ of them) being recomputed
by brute force from the definition of contiguity.  Rigidity passes to
larger subcomplexes by heredity.  On the hemi-icosahedron, where the two
classes are the two values of the $\Zt$ winding, each witness is also a
proof by Lemma~\ref{lem:winding}.  For example, for $D(1,1)$ and any facet
over $(\{0,4,5\},\{2,3,5\})$, the points $(0,0),(0,2),(5,5)$ span a hollow
triangle whose $\pi_1$-walk $0,0,5$ has winding $0$ and whose $\pi_2$-walk
$0,2,5$ runs around the essential empty triangle $\{0,2,5\}$.
\end{proof}

By the facet-class reduction of Section~\ref{sec:back}, maximality among
domains generated by facets is the notion that matters for covers.  The
same witness search succeeds for $16$ of the $49$ star rectangles of the
hemicube and for $4$ of the $64$ of the Klein bottle.

Theorem~\ref{thm:maximal} explains how the annealer gets stuck near the
rectangle family.  The rectangle chains $\pi_1\sim c_v\sim\pi_2$ have no
defect on the facets of their rectangle, so the search is drawn to them, and once a chain's domain
is a rectangle, no move can add a facet to it without first removing
others.  The hemi-icosahedron plateau at $k=\rho-1=5$ is of this kind.  Its
coverage $540/600$ is that of the union of the diagonal rectangles
$D(v,v)$, $v\ne5$; four of its domains are such rectangles ($150$ facets
each) and the fifth has $141$ facets.  The $60$ uncovered facets are those
covered only by $D(5,5)$, and none of them is an open door of any domain
(Table~\ref{tab:audit}).

\begin{table}[H]
\centering\footnotesize
\setlength{\tabcolsep}{3.5pt}
\begin{tabular}{@{}>{\raggedright\arraybackslash}p{4.3cm}crrr>{\raggedright\arraybackslash}p{2.2cm}>{\raggedright\arraybackslash}p{2.9cm}@{}}
\toprule
state & $k$ & covered & $|U|$ & dead & witnesses & exact completion \\
\midrule
\texttt{rp2\_k5\_hybrid} & $5$ & $540/600$ & $60$ & $60/60$ & hollow triangles & none ($\le12$ links)\\
\texttt{rp2\_k5\_from\_k6\_\dots7232} & $5$ & $544/600$ & $56$ & $54/56$ & hollow triangles & none\\
\texttt{rp2\_k3\_cold} & $3$ & $398/600$ & $202$ & $200/202$ & hollow triangles & ---\\
\texttt{hemicube\_k6\_v2} & $6$ & $836/864$ & $28$ & $27/28$ & hollow triangles & none\\
\texttt{klein\_k9\_v4} & $9$ & $1494/1536$ & $42$ & $41/42$ & $3$-walks, one $4$-walk & none\\
\texttt{klein\_k8\_r1} & $8$ & $1483/1536$ & $53$ & $53/53$ & $3$-, $4$-walks & ---\\
\texttt{hemi\_dodeca\_k6\_v2} & $6$ & $1922/1944$ & $22$ & $6/22$ & $3$-walks & none (open doors in five domains)\\
\texttt{mobius\_seed\_999} & $3$ & $148/150$ & $2$ & $0/2$ & --- & none\\
\texttt{mobius\_k3\_cold\_s4}, \texttt{\dots\_stretch\_from\_149} & $3$ & $149/150$ & $1$ & $0/1$ & --- & \textbf{found} ($8$ links)\\
\texttt{s3\_k2\_v2\_point}, \texttt{s3\_k2\_A}, \texttt{s3\_k2\_v2\_track} & $2$ & $460/500$ & $40$ & $40/40$ & rigid spheres & none ($\le8$ links)\\
\texttt{icosa\_k4\_v3} & $4$ & $2376/2400$ & $24$ & $0/24$ & --- & \textbf{found} ($8$ links): \texttt{icosa\_k4\_exact}\\
\bottomrule
\end{tabular}
\caption{Door audit of the stalled states (\texttt{rigid.py audit}; run
names are directories under \texttt{results/}).  ``Dead'' counts uncovered
facets that are an open door of no domain for the family of
Table~\ref{tab:floors}; a mod-$2$ cohomological balance test
(\texttt{scripts/relaxation/relax.py}) agrees on every row where it was run.  The last column is the exact
decision of Section~\ref{sec:exact} of whether $D_j\cup U$ is a domain for
some $j$, with all uncovered facets added at once.}
\label{tab:audit}
\end{table}

Table~\ref{tab:audit} separates the stalls into two kinds.  On the
hemi-icosahedron, the hemicube, the Klein bottle and $\partial\Delta^4$ the
states are \emph{rigid-closed}: all uncovered facets but at most two are
dead, so the state can improve only after some domain gives up facets.  On
the icosahedron and the M\"obius $149$-states every uncovered facet is an
open door, and the exact method of the next section completes the state.
On the hemi-dodecahedron and the M\"obius $148$-state the doors are open as
well, but no single domain absorbs all uncovered facets at once, and these
states remain unresolved.

\section{Exact methods}\label{sec:exact}

\paragraph{Deciding a single domain.}  Whether a facet set $D$ admits a
chain of at most $L$ links is a SAT instance: one-hot variables for the
values of the $L-1$ interior maps on the grid points of $D$, and, for every
facet and link, ``value present'' variables on which every minimal non-face
of $K$ is forbidden (for the complexes here these are non-edges, empty
triangles, and for $\partial\Delta^{n+1}$ the vertex set;
\texttt{scripts/rigid\_analysis/domain\_chain\_sat.py}).  Padding by $\pi_2$
makes ``no chain of exactly $L$ links'' equivalent to ``no chain of at most
$L$ links''.  Instances with up to $1\,000$ facets and $L=12$ solve in under
a second.  The encoder reproduces every recorded domain and every witness
of Theorem~\ref{thm:maximal} (UNSAT for every $L$ tested).

\paragraph{Completing a stall.}  For each domain $D_j$ of a stalled state
we ask whether $D_j\cup U$ is a domain.  On \texttt{icosa\_k4\_v3}, whose
four chain domains have $1099$, $881$, $744$ and $691$ facets (they
overlap), the answer is yes for the domain of $744$ facets at $L=8$, and no
for the smallest domain at $L\le12$.  Replacing that chain by the new one
gives a verified cover; after re-partitioning the facets its domains have
the sizes below.

\begin{proposition}\label{prop:icosa}
$\SCst(\text{icosahedron boundary})\le3$, by the certificate
\texttt{icosa\_k4\_exact} (in \texttt{results/}): four domains of $819$, $728$, $629$ and
$224$ facets with chains of $8$, $9$, $7$ and $6$ links.  Hence
$2\le\SCst(\text{icosahedron boundary})\le3$.
\end{proposition}

The same step completes both M\"obius $149$-states into $3$-covers (the one
domain with an open door absorbs the last facet with an $8$-link chain),
and it fails on every rigid-closed state of Table~\ref{tab:audit}, as it
must when every facet is dead.  The annealer's $k=4$ icosahedron plateau,
reached by three runs of $200$ epochs each, was therefore one chain away
from a cover; the annealer's local proposals did not find that chain.

\paragraph{A lazy exact cover solver.}  Combining the two halves gives a
solver for ``does $K\times K$ admit $k$ domains with chains of at most $L$
links'' (\texttt{scripts/rigid\_analysis/}\allowbreak\texttt{cegar\_cover.py}).
An outer SAT instance proposes $k$ disjoint facet classes avoiding the rigid
families; the single-domain decision checks each class; a class with no chain of at most $L$ links yields an unsatisfiable core,
which is minimised to a minimal facet set with no such chain and learned as
a clause for every class.  A SAT answer
returns chains and writes a certificate; an UNSAT answer is a statement for
chains of at most $L$ links, and the learned sets form a re-checkable
certificate.  The solver finds the three-domain M\"obius cover in one
second (the first proposal of the outer instance is already a cover), where
eight annealer runs stalled; it reports $k=2$ impossible on $K_M$ and on
$\partial\Delta^4$, where the rigid families alone are already
unsatisfiable; and it covers $\partial\Delta^3$ by three domains in two
seconds.  On the hemi-icosahedron at $k=5$, the hemicube at $k=5$ and the
icosahedron at $k=3$ it learns minimal sets of $3$ to $25$ facets with no
chain of at most $L$ links.  None of them contains a member of the rigid
families, since the proposed classes avoid those families.  Each iteration takes several
minutes, dominated by the outer instance; within a two-hour budget each run
reached $22$, $6$ and $17$ iterations with $290$, $90$ and $153$ learned
sets and no verdict, every proposed partition still containing a class
that is not a domain.  The learned sets are kept
(\texttt{certificates/learned\_*.json}) to seed later runs.  The solver is
effective at finding covers but weak at proving infeasibility, because its
learned cuts are local: on the cube boundary at $k=2$, where $\TC(S^2)=2$
already rules out a cover, $83$ iterations and $498$ learned sets did not
reach UNSAT in $38$ minutes.

\section{A direct CP-SAT model for \texorpdfstring{$\partial\Delta^4$}{the 3-sphere}}\label{sec:sat}

Before the rigid-sphere argument was available we tested
$\partial\Delta^4$ at $k=2$ with a direct CP-SAT model
(\texttt{paper-overleaf/paper/scs\_sat.py}): one vertex-valued variable per
grid vertex and interior map of each of $k$ chains, a facet covered by a
chain when every link's image union lies in a facet of $K$, and every facet
covered by some chain.  The encoding was validated on $1100$ random chains
against the annealer's defect counts, by round-tripping five recorded
covers, and by five solver controls in both directions
(\texttt{experiments/s3\_sat\_validation\_log.md}).  It returned INFEASIBLE
for $k=2$ at $L=4,5,8,12,16,24,32,40,48$ ($12$\,s to $1988$\,s,
\texttt{logs/sat\_ladder/}), a solver result for chains of at most $48$
links.  Theorem~\ref{thm:floors} now gives the statement for all lengths,
so the question whether every domain of
$\partial\Delta^4\times\partial\Delta^4$ admits a chain of at most $48$
links is no longer needed for $\SCst(\partial\Delta^4)=2$; every domain we
chained needed at most $4$ links.  The model does not scale: for
$\partial\Delta^6$ at $L=3,4$ and for the hemi-icosahedron at $k=5$, $L=4$
it exhausted $5$--$6$\,GB without a verdict, which is why the lazy solver of
Section~\ref{sec:exact} decides one domain at a time.

\section{Results}\label{sec:results}

\begin{table}[H]
\centering\footnotesize
\setlength{\tabcolsep}{3pt}
\begin{tabular}{@{}lrrrcc>{\raggedright\arraybackslash}p{4.3cm}c@{}}
\toprule
complex & $|V|$ & $|K{\times}K|$ & $\rho$ & $\TC$ & $\SCst\ge$ & best certified $\SCst\le$ (source) & audit \\
\midrule
circle $\partial\Delta^2$ & $3$ & $18$ & $3$ & $1$ & $1$ & $1$ (\texttt{circle\_probe3}); $=\TC$ & yes\\
wedge $S^1\vee S^1$ & $4$ & $50$ & $4$ & $2$\rlap{$^*$} & $2$ & $2$ (\texttt{wedge\_k3\_v2}); $=\TC$ & yes\\
$\partial\Delta^{n+1}$, $n$ even & & & $3$ & $2$ & $2$ & $2$; $=\TC$ (Cor.~\ref{cor:rho}) & yes ($n\le8$)\\
$\partial\Delta^{n+1}$, $n=3,5,7,9$ & & & $3$ & $1$ & $\mathbf 2$ & $2$ (Cor.~\ref{cor:rho}); $=\SCst>\TC$ (Thm.~\ref{thm:floors}) & yes ($n\le7$)\\
M\"obius band & $5$ & $150$ & $4$ & $1$ & $\mathbf 2$ & $2$ (\texttt{mobius\_k3\_sat}, \texttt{mobius\_k3\_cegar\_L8}); $=\SCst$ & yes\\
hemi-icosahedron & $6$ & $600$ & $6$ & $3$ & $\mathbf 4$ & $5$ (\texttt{rp2\_k6}) & yes\\
hemicube & $7$ & $864$ & $7$ & $3$ & $\mathbf 4$ & $6$ (\texttt{hemicube\_k7}) & internal\\
hemi-octahedron & $9$ & $1536$ & $9$ & $3$ & $3$ & $7$ (\texttt{hemi\_octa\_k8}) & internal\\
hemi-dodecahedron & $10$ & $1944$ & $15$ & $3$ & $3$ & $6$ (\texttt{hemi\_dodeca\_k7}) & internal\\
Klein bottle & $8$ & $1536$ & $12$ & $4$ & $4$ & $9$ (\texttt{klein\_k10}) & yes\\
Cs\'asz\'ar torus & $7$ & $1176$ & $11$ & $2$ & $\mathbf 3$ & $10$ (star cover) & yes\\
cube boundary & $8$ & $864$ & $8$ & $2$ & $2$ & $2$ (\texttt{cube\_k3\_guided}); $=\TC$ & yes\\
icosahedron boundary & $12$ & $2400$ & $24$ & $2$ & $2$ & $3$ (\texttt{icosa\_k4\_exact}) & yes\\
dodecahedron boundary & $20$ & $7776$ & $36$ & $2$ & $2$ & $35$ (star cover) & yes\\
\bottomrule
\end{tabular}
\caption{State of knowledge at level $0$.  The column ``$\SCst\ge$'' is the
best lower bound on $\SCst$ (one less than the number of domains in
Table~\ref{tab:floors}), bold where it exceeds $\TC$
(Theorem~\ref{thm:floors}); $\TC$ is reduced \cite{F03,FTY03,CV17}; $^*$for the wedge see
\cite[Thm.~7.3]{F04}.
``Audit'': the certificate was re-checked by \texttt{src/independent\_check.py};
for the simplex boundaries the archived star covers go up to
$\partial\Delta^9$, and Corollary~\ref{cor:rho} is a proof in any case.
``Internal'': the auditor cannot yet rebuild the complex (its facet counts
$864$, $1536$ and $1944$ collide with other complexes in, or are absent
from, its table), so only the search's own verifier was run.}
\label{tab:surfaces}
\end{table}

Four observations.  First, $\SCst=\TC$ on the circle, the wedge, all even
simplex boundaries and the cube boundary, and $\SCst=\TC+1$ on
$\partial\Delta^{n+1}$ for $n=3,5,7,9$ and on the M\"obius band; this gap at
level $0$ holds for chains of every length.  Second, the four
triangulations of $S^2$ behave differently: $\partial\Delta^3$ and the cube
have $\SCst=2$, the icosahedron has $2\le\SCst\le3$, the dodecahedron only
$\SCst\le35$, while $\rho$ is $3$, $8$, $24$ and $36$.  Third, the two-chain
plateaus $460$ on $\partial\Delta^4$ and $2380$ on $\partial\Delta^5$ equal
the number of facets covered by the best pair of star rectangles
($23\cdot20$ and $34\cdot70$), and Theorem~\ref{thm:maximal} explains why
rectangle-shaped domains stop growing.  Fourth, on the hemi-octahedron and
the hemi-dodecahedron the rigid families give nothing beyond $\TC+1=4$
domains, while the best covers have $8$ and $7$; on the torus the families
give $4$ domains and no cover with fewer than $\rho=11$ is known.

\section{Certificates, auditing, reproduction}\label{sec:repro}

Every cover reported is stored as data: chains are replayed
deterministically, reduced by \cite[Alg.~2]{OLGB21}, and checked to connect
the two projections through pairwise contiguous simplicial maps over every
facet of every domain (\texttt{cover.txt}, \texttt{planners.txt},
\texttt{summary.json}).  An independent auditor
(\texttt{src/independent\_check.py}), a pure-Python program that rebuilds
each triangulation from its facet list and shares no code with the search,
repeats every check and exits with an error on any discrepancy.  It
verifies $59$ archived covers: the star covers of every complex of
Table~\ref{tab:rho} except the wedge and the three subdivided models of
$\mathbb{RP}^2$ (for the simplex boundaries up to $\partial\Delta^9$, with
$1\,287\,000$ facets), the whole icosahedron ladder (the $24$-domain seed,
$k=23,\dots,5$ and \texttt{icosa\_k4\_exact}), the hemi-icosahedron and
Klein rungs, the M\"obius certificates and the certificates of
Table~\ref{tab:annealer} (\texttt{audit/verify\_all.py --check} in the
certificates repository, which also checks a hash manifest).  The hemicube,
hemi-octahedron and hemi-dodecahedron covers pass only the internal verifier
(Table~\ref{tab:surfaces}); they are archived and hashed separately.

The computer-assisted theorems rest on three kinds of artefact: SAT
instances built from proved rigid families (\texttt{rigid.py bound}, logs
in \texttt{scripts/rigid\_analysis/}\allowbreak\texttt{certificates/}\allowbreak\texttt{bound\_*.txt},
and the rigid-sphere scripts in \texttt{scripts/relaxation/}); witness
lists re-verified from the definitions
(\texttt{maximal\_rectangles\_*.json}); and single-domain decisions, each
re-checkable in seconds.  None of the solver runs produced checked proof
logs.  The CP-SAT infeasibilities of Section~\ref{sec:sat} are solver
outputs and are labelled as such.

Source code (annealer, star cover, rigid-subcomplex tools, exact solvers,
auditor) is at the address below; file paths in this paper are relative to
its directory \texttt{opti/} unless they begin with \texttt{opti/}.
\begin{center}\url{https://github.com/MiguelPerezI/chain-first-annealer}\end{center}
and the archived certificates at
\begin{center}\url{https://github.com/MiguelPerezI/simplicial-complexity-certificates}\end{center}
Every cover cited in this paper is in the certificates archive; the
rigid-analysis certificates, logs and stalled states are in
\texttt{opti/scripts/rigid\_analysis/}, \texttt{opti/results/} and
\texttt{opti/experiments/} of the source repository.  The versions described
here carry the tag \texttt{arxiv-v1} in both repositories.  All annealer
runs except $S^6$ used one Tesla T4; wall times are single-run
observations.  The rigid-analysis and exact computations ran on a
four-core laptop CPU.

\section{Concluding remarks}\label{sec:concl}

The annealer has two specific weaknesses.  It does not find completions
that require re-routing a whole chain at once (the icosahedron and the
M\"obius $149$-states, completed by the exact method), and it is drawn into
the rectangle family, whose members are maximal domains on three of the
surfaces.  Where a stalled state is rigid-closed, no extension of its
domains by single facets exists, and the remaining gap is between the
certified covers and the rigid-family lower bounds of
Table~\ref{tab:surfaces}.

All complexes in this paper are minimal, so on all of them every graph
$\Gamma_w$ is rigid.  What distinguishes the odd spheres, the M\"obius band,
the torus and two of the triangulations of $\mathbb{RP}^2$ is that their
small rigid pieces cannot be avoided by $\TC+1$ classes of facets; on the
icosahedron no short closed walk is rigid and the graphs can be avoided by
two classes.  Every rigid graph $\Gamma_w$ we tested on these complexes
joins non-homotopic maps (Remark~\ref{rem:graphs}(c)), but rigid closed
walks need not: on the Klein bottle the hollow triangle
$(0,1),(2,3),(5,6)$ is rigid because its projected $3$-walks $0,2,5$ and
$1,3,6$ lie in different contiguity classes, yet after repeating the first
vertex the two $4$-walks lie in the same class, so the projections are
homotopic.  Whether the level-$0$ excess on the odd spheres and the
surfaces is detected by homotopy-theoretic rigidity alone is open.

Open questions: whether the hemi-icosahedron admits five domains
($\SCst\in\{4,5\}$), whether the hemicube admits five or six
($\SCst\in\{4,5,6\}$), whether the icosahedron admits three
($\SCst\in\{2,3\}$); whether the Klein bottle needs more than five domains
(the rigid family at $k=5$ is undecided); the torus, which needs between
$4$ and $11$ domains; and whether $\SCst(\partial\Delta^{n+1})=2$ for every
odd $n\ge3$, which is known here for $n\le9$.

\end{document}